\documentclass{amsart}
\usepackage{amsmath,amssymb}
\usepackage{yhmath}

\newcommand{\bdis}{\begin{displaymath}}
\newcommand{\edis}{\end{displaymath}}
\newcommand{\be}{\begin{equation}}
\newcommand{\ee}{\end{equation}}
\newcommand{\mbb}{\mathbb}
\newcommand{\mcal}{\mathcal}

\newcommand{\vp}{\varphi} 
\newcommand{\vth}{\vartheta}

\newcommand{\zf}{\zeta\left(\frac{1}{2}+it\right)}

\newcommand{\FR}{\frac{x^n+y^n}{z^n}} 
 
\newcommand{\FRm}{\frac{x^m+y^m}{z^m}}

\DeclareMathOperator{\im}{Im}

\theoremstyle{definition}

\theoremstyle{remark}
\newtheorem{remark}[]{Remark}

\newtheorem*{mydef11}{{\bf Theorem 1}}

\newtheorem*{mydef12}{{\bf Theorem 2}}

\newtheorem*{mydef13}{{\bf Theorem 3}}

\newtheorem*{mydef14}{{\bf Theorem 4}}

\newtheorem*{mydef15}{{\bf Theorem 5}} 

\newtheorem*{mydef16}{{\bf Theorem 6}} 

\newtheorem*{mydef17}{{\bf Theorem 7}}

\newtheorem*{mydef41}{{\bf Corollary 1}}

\newtheorem*{mydef42}{{\bf Corollary 2}}

\newtheorem*{mydef51}{{\bf Lemma 1}}

\newtheorem*{mydef52}{{\bf Lemma 2}}

\newtheorem*{mydef53}{{\bf Lemma 3}}

\newtheorem*{mydef81}{{\bf Property 1}}

\newtheorem*{mydef82}{{\bf Property 2}}

\newtheorem*{mydef83}{{\bf Property 3}} 

\newtheorem*{mydef84}{{\bf Property 4}} 

\newtheorem*{mydef85}{{\bf Property 5}} 

\newtheorem*{mydef86}{{\bf Property 6}} 

\newtheorem*{mydef87}{{\bf Property 7}}

\numberwithin{equation}{section}

\begin{document}

\title[Jacob's ladders and new $\zeta$-functionals \dots]{Jacob's ladders and new $\zeta$-functionals and corresponding $\zeta$-equivalents of the Fermat-Wiles theory based on sums of elementary $\zeta$-pulses}

\author{Jan Moser}

\address{Department of Mathematical Analysis and Numerical Mathematics, Comenius University, Mlynska Dolina M105, 842 48 Bratislava, SLOVAKIA}

\email{jan.mozer@fmph.uniba.sk}


\begin{abstract}
In this paper we obtain some sum-variants of our first $\zeta$-equivalent of the Fermat-Wiles theorem. For example, on the Riemann hypothesis, we give a new type of $\zeta$-equivalent based on a sum of isolated rectangle-shaped signals (pulses).  
\end{abstract}
\maketitle

\section{Introduction} 

\subsection{} 

Let us remind the classical Hardy-Littlewood integral (1918)\footnote{See \cite{1}.} 
\be \label{1.1} 
J(T)=\int_0^T\left|\zf\right|^2{\rm d}t=\int_0^T Z^2(t){\rm d}t, 
\ee 
where we use also Riemann's real function 
\bdis 
\begin{split}
& Z(t)=e^{i\vth(t)}\zf \ \Rightarrow \ |Z(t)|=\left|\zf\right|, \\ 
& \vth(t)=-\frac{t}{2}\ln\pi+\im\ln\Gamma\left(\frac 14+i\frac t2\right). 
\end{split}
\edis  
Next, 105 years after, we have obtained a new class of increments of the Hardy-Littlewood integral (\ref{1.1}). Namely, we have proved the existence of the almost linear increments 
\be \label{1.2} 
\int_T^{\overset{1}{T}(T)}Z^2(t){\rm d}t\sim (1-c)T,\ T\to\infty, 
\ee  
where 
\be \label{1.3} 
\overset{1}{T}(T)=\vp_1^{-1}(T), 
\ee 
see \cite{7}, (3.4), (3.6), $r=1$, and $c$ is the Euler's constant. 

\subsection{} 

Consequently, we have constructed on the basis of formula (\ref{1.2}) the following $\zeta$-functional\footnote{See \cite{8}, (4.6) -- (4.9).} 
\be \label{1.4} 
\lim_{\tau\to\infty}\frac{1}{\tau}\int_{\frac{x}{1-c}\tau}^{[\frac{x}{1-c}\tau]^1}Z^2(t){\rm d}t=x
\ee 
for every fixed $x>0$. To be precise, we call the mapping $\mcal{F}_1$: 
\be \label{1.5} 
y(\tau,x)=\frac{x}{1-c}\tau\xrightarrow{\mcal{F}_1}x
\ee 
defined by (\ref{1.4}) as $\zeta$-functional. 

In the special case of Fermat's rationals 
\be \label{1.6} 
x\to \FR,\ x,y,z,n\in\mbb{N},\ n\geq 3
\ee 
one obtains the following formula 
\be \label{1.7} 
\begin{split} 
& \lim_{\tau\to\infty}\frac{1}{\tau}\int_{\FR\frac{\tau}{1-c}}^{[\FR\frac{\tau}{1-c}]^1}Z^2(t){\rm d}t=\FR,\\ 
& [G]^1=\vp_1^{-1}(G), 
\end{split} 
\ee 
see \cite{8}, (5.2). 

Consequently, we have obtained this result\footnote{See \cite{8}, (5.3).}: The $\zeta$-condition  
\be \label{1.8} 
\lim_{\tau\to\infty}\frac{1}{\tau}\int_{\FR\frac{\tau}{1-c}}^{[\FR\frac{\tau}{1-c}]^1}Z^2(t){\rm d}t\not=1
\ee 
on the class of all Fermat's rationals represents the $\zeta$-equivalent of the Fermat-Wiles theorem. 

\begin{remark}
Formula (\ref{1.8}) has been the first point of contact between our almost linear formula (\ref{1.2}) and the Fermat-Wiles theorem, i.e. between Riemann zeta-function and the Fermat-Wiles theorem. The series of our sixteen papers issued subsequently after the work \cite{8} contain set of further points of contact with the Fermat-Wiles theorem. 
\end{remark} 

\subsection{} 

The integral in the formula (\ref{1.4}) concealed his internal structure. Our aim is to obtain more \emph{stable} (or coded) form - namely a sum form of this integral. 

Let the symbol $\mcal{A}_0(n)$ denote the region bounded by the graph of the function 
\be \label{1.9} 
Z^2(t),\ t\in [\gamma_n,\gamma_{n+1}];\ Z(\gamma_n)=Z(\gamma_{n+1})=0, 
\ee   
and the segment $[\gamma_n,\gamma_{n+1}]$.\footnote{The symbol $\{\gamma\}$ stands for the set of zeros of the function $Z(t)$.} For the area $|\mcal{A}_0(n)|$ of this region we have, of course, 
\be \label{1.10} 
|\mcal{A}_0(n)|=\int_{\gamma_n}^{\gamma_{n+1}}Z^2(t){\rm d}t. 
\ee  
In this case we get the following functional 
\be \label{1.11} 
\lim_{\tau\to\infty}\frac{1}{\tau}\left\{
\sum_{\gamma_n>\frac{x}{1-c}\tau}^{\gamma_n<[\frac{x}{1-c}\tau]^1}\int_{\gamma_n}^{\gamma_{n+1}}Z^2(t){\rm d}t
\right\}=x 
\ee 
for every fixed $x>0$. 

\begin{remark}
The formula (\ref{1.11}) represents the first sum-variant of the basic formula (\ref{1.4}) for the corresponding set 
\be \label{1.12} 
\bigcup_{n=1}^N \mcal{A}_0(n)
\ee 
of the $\zeta$-curvilinear signals $\mcal{A}_0(n)$. 
\end{remark} 

In the special case of the Fermat's rationals 
\be \label{1.13} 
x\to \FRm,\ x,y,z,m\in\mbb{N},\ m\geq 3
\ee 
we obtain the following result: the $\zeta$-condition 
\be \label{1.14} 
\lim_{\tau\to\infty}\frac{1}{\tau}\left\{
\sum_{\gamma_n>\FRm\frac{\tau}{1-c}}^{\gamma_n<[\FRm\frac{\tau}{1-c}]^1}\int_{\gamma_n}^{\gamma_{n+1}}Z^2(t){\rm d}t
\right\}\not=1
\ee 
on the set of all Fermat's rationals expresses the new $\zeta$-equivalent of the Fermat-Wiles theorem. 

\subsection{} 

Let 
\be \label{1.15} 
t_0(n):\ Z'(t_0(n))=0,\ \gamma_n<t_0(n)<\gamma_{n+1}, 
\ee  
and next, let values $\xi_1(m)$ and $\xi_2(m)$, where 
\be \label{1.16} 
\begin{split}
& \xi_1=t_0-\frac{1}{Z^2(t_0)}\int_{\gamma_n}^{t_0}Z^2(t){\rm d}t, \\ 
& \xi_2=t_0+\frac{1}{Z^2(t_0)}\int_{t_0}^{\gamma_{n+1}}Z^2(t){\rm d}t, 
\end{split} 
\ee 
represent the Bonnet points of the integral 
\be \label{1.17} 
\int_{\gamma_n}^{\gamma_{n+1}}Z^2(t){\rm d}t. 
\ee  
Then, the influence of the Riemann hypothesis on the functional behaves as follows: on Riemann hypothesis we have the following $\zeta$-functional in the sum-form 
\be \label{1.18} 
\lim_{\tau\to\infty}\frac{1}{\tau}\left\{
\sum_{t_0(n)>\frac{x}{1-c}\tau}^{t_0(n)<[\frac{x}{1-c}\tau]^1}Z^2(t_0(n))[\xi_2(n)-\xi_1(n)]
\right\}=x 
\ee 
for every fixed $x>0$, where 
\be \label{1.19} 
\gamma_n<\xi_1(n)<t_0(n)<\xi_2(n)<\gamma_{n+1},\ n=1,\dots,N. 
\ee 

Next, in the special case of the Fermat's rationals, see (\ref{1.13}), we have the following result. On Riemann hypothesis it is true that the $\zeta$-condition: 
\be \label{1.20} 
\begin{split}
& \lim_{\tau\to\infty}\frac{1}{\tau}\left\{
\sum_{t_0(n)>\FRm\frac{\tau}{1-c}}^{t_0(n)<[\FRm\frac{\tau}{1-c}]^1}Z^2(t_0(n))[\xi_2(n)-\xi_1(n)]
\right\}\not=1 
\end{split}
\ee 
on the set of all Fermat's rationals expresses new $\zeta$-equivalent of the Fermat-Wiles theorem. 

Let the symbol $\mcal{B}(n)$ stands for the rectangle 
\be \label{1.21} 
(\xi_1(n),\xi_2(n))\times (0,Z^2(t_0)), 
\ee  
where, of course, 
\be \label{1.22} 
|\mcal{B}(n)|=Z^2(t_0)[\xi_2(n)-\xi_1(n)]. 
\ee 

\begin{remark}
Consequently, by the Riemann hypothesis, we have alternated the set\footnote{See (\ref{1.12}).}
\bdis 
\bigcup_{n=1}^N \mcal{A}_0(n)
\edis 
of curvilinear signals $\mcal{A}_0(n)$ by the disconnected set\footnote{See (\ref{1.19}) and (\ref{1.21}).} 
\be \label{1.23} 
\bigcup_{n=1}^N\mcal{B}(n) 
\ee  
of the isolated rectangular pulses, see (\ref{1.18}). 
\end{remark} 

\begin{remark}
Let us remind that for the construction of the $\zeta$-equivalent (\ref{1.20}) of the Fermat-Wiles theorem on the set of $\zeta$-pulses $\mcal{B}(n)$ we have used the same consequences of the Riemann hypothesis as we have used 52 years ago in our paper \cite{2}, see also \cite{3}, and \cite{9}, pp. 127-132, for the construction of the first infinite set of mathematical models of Universe in the relativistic cosmology based on the Riemann zeta-function. 
\end{remark}

\section{Jacob's ladders: notions and basic geometrical properties}  

\subsection{}

In this paper we use the following notions of our works \cite{4} -- \cite{8}: 
\begin{itemize}
\item[{\tt (a)}] Jacob's ladder $\vp_1(T)$, 
\item[{\tt (b)}] direct iterations of Jacob's ladders 
\bdis 
\begin{split}
	& \vp_1^0(t)=t,\ \vp_1^1(t)=\vp_1(t),\ \vp_1^2(t)=\vp_1(\vp_1(t)),\dots , \\ 
	& \vp_1^k(t)=\vp_1(\vp_1^{k-1}(t))
\end{split}
\edis 
for every fixed natural number $k$, 
\item[{\tt (c)}] reverse iterations of Jacob's ladders 
\be \label{2.1}  
\begin{split}
	& \vp_1^{-1}(T)=\overset{1}{T},\ \vp_1^{-2}(T)=\vp_1^{-1}(\overset{1}{T})=\overset{2}{T},\dots, \\ 
	& \vp_1^{-r}(T)=\vp_1^{-1}(\overset{r-1}{T})=\overset{r}{T},\ r=1,\dots,k, 
\end{split} 
\ee   
where, for example, 
\be \label{2.2} 
\vp_1(\overset{r}{T})=\overset{r-1}{T}
\ee  
for every fixed $k\in\mbb{N}$ and every sufficiently big $T>0$. We also use the properties of the reverse iterations listed below.  
\be \label{2.3}
\overset{r}{T}-\overset{r-1}{T}\sim(1-c)\pi(\overset{r}{T});\ \pi(\overset{r}{T})\sim\frac{\overset{r}{T}}{\ln \overset{r}{T}},\ r=1,\dots,k,\ T\to\infty,  
\ee 
\be \label{2.4} 
\overset{0}{T}=T<\overset{1}{T}(T)<\overset{2}{T}(T)<\dots<\overset{k}{T}(T), 
\ee 
and 
\be \label{2.5} 
T\sim \overset{1}{T}\sim \overset{2}{T}\sim \dots\sim \overset{k}{T},\ T\to\infty.   
\ee  
\end{itemize} 

\begin{remark}
	The asymptotic behaviour of the points 
	\bdis 
	\{T,\overset{1}{T},\dots,\overset{k}{T}\}
	\edis  
	is as follows: at $T\to\infty$ these points recede unboundedly each from other and all together are receding to infinity. Hence, the set of these points behaves at $T\to\infty$ as one-dimensional Friedmann-Hubble expanding Universe. 
\end{remark}  

\subsection{} 

Let us remind that we have proved\footnote{See \cite{7}, (3.4).} the existence of almost linear increments 
\be \label{2.6} 
\begin{split}
& \int_{\overset{r-1}{T}}^{\overset{r}{T}}\left|\zf\right|^2{\rm d}t\sim (1-c)\overset{r-1}{T}, \\ 
& r=1,\dots,k,\ T\to\infty,\ \overset{r}{T}=\overset{r}{T}(T)=\vp_1^{-r}(T)
\end{split} 
\ee 
for the Hardy-Littlewood integral (1918), \cite{1}: 
\be \label{2.7} 
J(T)=\int_0^T\left|\zf\right|^2{\rm d}t. 
\ee  

For completeness, we give here some basic geometrical properties related to Jacob's ladders. These are generated by the sequence 
\be \label{2.8} 
T\to \left\{\overset{r}{T}(T)\right\}_{r=1}^k
\ee 
of reverse iterations of the Jacob's ladders for every sufficiently big $T>0$ and every fixed $k\in\mbb{N}$. 

\begin{mydef81}
The sequence (\ref{2.8}) defines a partition of the segment $[T,\overset{k}{T}]$ as follows 
\be \label{2.9} 
|[T,\overset{k}{T}]|=\sum_{r=1}^k|[\overset{r-1}{T},\overset{r}{T}]|
\ee 
on the asymptotically equidistant parts 
\be \label{2.10} 
\begin{split}
& \overset{r}{T}-\overset{r-1}{T}\sim \overset{r+1}{T}-\overset{r}{T}, \\ 
& r=1,\dots,k-1,\ T\to\infty. 
\end{split}
\ee 
\end{mydef81} 

\begin{mydef82}
Simultaneously with the Property 1, the sequence (\ref{2.8}) defines the partition of the integral 
\be \label{2.11} 
\int_T^{\overset{k}{T}}\left|\zf\right|^2{\rm d}t
\ee 
into the parts 
\be \label{2.12} 
\int_T^{\overset{k}{T}}\left|\zf\right|^2{\rm d}t=\sum_{r=1}^k\int_{\overset{r-1}{T}}^{\overset{r}{T}}\left|\zf\right|^2{\rm d}t, 
\ee 
that are asymptotically equal 
\be \label{2.13} 
\int_{\overset{r-1}{T}}^{\overset{r}{T}}\left|\zf\right|^2{\rm d}t\sim \int_{\overset{r}{T}}^{\overset{r+1}{T}}\left|\zf\right|^2{\rm d}t,\ T\to\infty. 
\ee 
\end{mydef82} 

It is clear, that (\ref{2.10}) follows from (\ref{2.3}) and (\ref{2.5}) since 
\be \label{2.14} 
\overset{r}{T}-\overset{r-1}{T}\sim (1-c)\frac{\overset{r}{T}}{\ln \overset{r}{T}}\sim (1-c)\frac{T}{\ln T},\ r=1,\dots,k, 
\ee  
while our eq. (\ref{2.13}) follows from (\ref{2.6}) and (\ref{2.5}).

\section{A $\zeta$-sum variants of almost linear formula (2023)} 

\subsection{} 

Let us remind our almost linear formula, \cite{7}, (3.4), (3.6), 
\be \label{3.1} 
\int_T^{\overset{1}{T}}\left|\zf\right|^2{\rm d}t=(1-c)T+\mcal{O}(T^{1/3+\delta}),\ T\to\infty, 
\ee  
where $\delta>0$ is any fixed small value. Next, let the symbol 
\be \label{3.2} 
\{\gamma_k\}_{k=1}^\infty, \ 0<\gamma_k<+\infty
\ee 
denote the set of all zeros of the function 
\be \label{3.3} 
\zf,\ 0<t<+\infty . 
\ee  
Now, let the symbol 
\be \label{3.4} 
\{\gamma_n\}_{n=1}^{N+1} 
\ee 
denote the subset of those zeros that is selected from the set (\ref{3.2}) by the interval 
\be \label{3.5} 
(T,\overset{1}{T}(T))
\ee 
as it follows 
\be \label{3.6} 
\begin{split}
& \gamma_0\leq T<\gamma_1<\gamma_2<\dots<\gamma_N<\gamma_{N+1}<\overset{1}{T}(T)\leq \gamma_{N+2}, \\ 
& N=N(T). 
\end{split}
\ee 

\begin{remark}
Of course, the letter $n$ stands for the relative index, comp. (\ref{3.2}). 
\end{remark} 

\begin{remark}
Let us remind that the part (\ref{3.2}) -- (\ref{3.6}) contains definition of the following set-valued function $F$: 
\be \label{3.7} 
(T,\overset{1}{T}(T))\xrightarrow{F}\{\gamma_n\}_{n=1}^{N+1},\ N=N(T),\ T\to\infty, 
\ee 
that is defined on the set 
\be \label{3.8} 
\{(T,\overset{1}{T}(T))\},\ T\to\infty 
\ee 
of the intervals, i.e. the finite set $\{\gamma_n\}_{n=1}^{N+1}$ of zeros with the property (\ref{3.6}) corresponds to every interval $(T,\overset{1}{T}(T))$. 
\end{remark} 

\subsection{} 

In the next we use also the equality, see Introduction, 
\bdis 
\left|\zf\right|=|Z(t)|. 
\edis  
We obtain, see (\ref{3.1}) and (\ref{3.6}), 
\be \label{3.9} 
\begin{split}
& \int_T^{\overset{1}{T}(T)}\left|\zf\right|^2{\rm d}t=\sum_{n=1}^N\int_{\gamma_n}^{\gamma_{n+1}}Z^2(t){\rm d}t+ \\ 
& \int_T^{\gamma_1}Z^2(t){\rm d}t+\int_{\gamma_{N+1}}^{\overset{1}{T}(T)}Z^2(t){\rm d}t. 
\end{split}
\ee 
For our purposes it is sufficient to use elementary estimates  
\be \label{3.10} 
\begin{split}
& \gamma_1-T,\ \overset{1}{T}(T)-\gamma_{N+1}=\mcal{O}(T^{1/6}), \\ 
& Z^2(t)=\mcal{O}(T^{1/3}),\ t\in[\gamma_0,\gamma_1]\cup [\gamma_{N+1},\gamma_{N+2}]. 
\end{split}
\ee 
These estimates imply the following\footnote{See (\ref{3.9}).} asymptotic formula 
\be \label{3.11} 
\begin{split}
& \int_T^{\overset{1}{T}(T)}\left|\zf\right|^2{\rm d}t=\sum_{n=1}^N\int_{\gamma_n}^{\gamma_{n+1}}Z^2(t){\rm d}t+\mcal{O}(\sqrt{T})= \\ 
& \sum_{T<\gamma_n<\overset{1}{T}(T)}\int_{\gamma_n}^{\gamma_{n+1}}Z^2(t){\rm d}t+\mcal{O}(\sqrt{T}),\ T\to\infty, 
\end{split}
\ee 
where, of course, $n=1,\dots,N$, comp. (\ref{3.6}). 

Consequently, by (\ref{3.1}) and (\ref{3.11}), the following lemma holds true. 

\begin{mydef51}
\be \label{3.12} 
\sum_{T<\gamma_n<\overset{1}{T}(T)}\int_{\gamma_n}^{\gamma_{n+1}}Z^2(t){\rm d}t=(1-c)T+\mcal{O}(\sqrt{T}),\ T\to\infty, 
\ee 
where, of course, $n=1,\dots,N(T)$. 
\end{mydef51} 

\begin{remark}
The formula (\ref{3.12}) represents the first sum-variant of the basic formula (\ref{3.1}) which is formulated for the corresponding set 
\be \label{3.13} 
\bigcup_{n=1}^N\mcal{A}_0(n),\ N=N(T) 
\ee 
of the curvilinear signals $\mcal{A}_0(n)$. 
\end{remark} 

\subsection{} 

Next, we use the usual mean value theorem in (\ref{3.12}) that gives us the following 
\be \label{3.14} 
\int_{\gamma_n}^{\gamma_{n+1}}Z^2(t){\rm d}t=Z^2(\bar{t})(\gamma_{n+1}-\gamma_n),\ \bar{t}\in(\gamma_n,\gamma_{n+1}), 
\ee  
where 
\be \label{3.15} 
\bar{t}=\bar{t}(n)=\bar{t}(\gamma_n,\gamma_{n+1};Z^2). 
\ee 
Now, it follows from (\ref{3.12}) that the next lemma holds true. 

\begin{mydef52} 
\be \label{3.16} 
\begin{split}
& \sum_{T<\gamma_n<\overset{1}{T}(T)}Z^2(\bar{t}(n))(\gamma_{n+1}-\gamma_n)= (1-c)T+\mcal{O}(\sqrt{T}),\ T\to\infty. 
\end{split}
\ee 
\end{mydef52} 

Now, let the symbol $\mcal{A}_1(n)$ denote the rectangle 
\be \label{3.17} 
[\gamma_n,\gamma_{n+1}]\times [0,Z^2(\bar{t}(n))], 
\ee  
where\footnote{See (\ref{1.10}) and (\ref{3.14}).} 
\be \label{3.18} 
|\mcal{A}_1(n)|=Z^2(\bar{t}(n))(\gamma_{n+1}-\gamma_n)=|\mcal{A}_0(n)|. 
\ee  

\begin{remark}
The formula (\ref{3.16}) represents the second sum-variant of our basic formula (\ref{3.1}). Namely, the formula (\ref{3.16}) contains the sum 
\be \label{3.19} 
\sum_{n=1}^N|\mcal{A}_1(n)| 
\ee  
of the areas of elementary rectangular signals $\mcal{A}_1(n)$. 
\end{remark}

\section{New sum-variants of our first $\zeta$-functional (2023) and corresponding $\zeta$-equivalents of the Fermat-Wiles theorem} 

\subsection{} 

If we use the substitution 
\be \label{4.1} 
T=\frac{x}{1-c}\tau,\ \{T\to+\infty\}\Leftrightarrow \{\tau\to+\infty\}
\ee  
for every fixed $x>0$ in the formula (\ref{3.12}), then we obtain the following functional. 

\begin{mydef11}
\be \label{4.2} 
\lim_{\tau\to\infty}\frac{1}{\tau}\left\{
\sum_{\gamma_n>\frac{x}{1-c}\tau}^{\gamma_n<[\frac{x}{1-c}\tau]^1}\int_{\gamma_n}^{\gamma_{n+1}}Z^2(t){\rm d}t
\right\}=x
\ee 
for fixed every $x>0$. 
\end{mydef11} 

\begin{remark}
The formula (\ref{4.2}) represents the first sum-variant of our original $\zeta$-functional\footnote{See \cite{8}, (4.6).}
\be \label{4.3} 
\lim_{\tau\to\infty}\frac{1}{\tau}
\int_{\frac{x}{1-c}\tau}^{[\frac{x}{1-c}\tau]^1}Z^2(t){\rm d}t
=x,\ x>0. 
\ee 
The new sum-variant (\ref{4.2}) expresses more stable analytical form than (\ref{4.3}). For example the integral mean-value theorem does not apply on the sum in (\ref{4.2}). 
\end{remark} 

\subsection{} 

In the special case of Fermat's rationals, see (\ref{1.13}), we obtain from (\ref{4.2}) this result. 

\begin{mydef41}
\be \label{4.4} 
\lim_{\tau\to\infty}\frac{1}{\tau}\left\{
\sum_{\gamma_n>\FRm\frac{\tau}{1-c}}^{\gamma_n<[\FRm\frac{\tau}{1-c}]^1}\int_{\gamma_n}^{\gamma_{n+1}}Z^2(t){\rm d}t
\right\}=\FRm
\ee 
for every fixed Fermat's rational. 
\end{mydef41} 

Consequently, the following theorem holds true. 

\begin{mydef12}
The $\zeta$-condition 
\be \label{4.5} 
\lim_{\tau\to\infty}\frac{1}{\tau}\left\{
\sum_{\gamma_n>\FRm\frac{\tau}{1-c}}^{\gamma_n<[\FRm\frac{\tau}{1-c}]^1}\int_{\gamma_n}^{\gamma_{n+1}}Z^2(t){\rm d}t
\right\}\not=1
\ee 
on the class of all Fermat's rationals expresses the new $\zeta$-equivalent of the Fermat-Wiles theorem. 
\end{mydef12} 

\subsection{} 

Here we give a list of similar results following from our Lemma 2. 

\begin{mydef13}
\be \label{4.6} 
\lim_{\tau\to\infty}\frac{1}{\tau}\left\{
\sum_{\gamma_n>\frac{x}{1-c}\tau}^{\gamma_n<[\frac{x}{1-c}\tau]^1}Z^2(\bar{t}(n))(\gamma_{n+1}-\gamma_n)
\right\}=x
\ee 
for every fixed $x>0$. 
\end{mydef13} 

\begin{mydef42}
\be \label{4.7} 
\lim_{\tau\to\infty}\frac{1}{\tau}\left\{
\sum_{\gamma_n>\FRm\frac{\tau}{1-c}}^{\gamma_n<[\FRm\frac{\tau}{1-c}]^1}Z^2(\bar{t}(n))(\gamma_{n+1}-\gamma_n)
\right\}=\FRm
\ee 
for every fixed Fermat's rational. 
\end{mydef42} 

\begin{mydef14}
The $\zeta$-condition 
\be \label{4.8} 
\lim_{\tau\to\infty}\frac{1}{\tau}\left\{
\sum_{\gamma_n>\FRm\frac{\tau}{1-c}}^{\gamma_n<[\FRm\frac{\tau}{1-c}]^1}Z^2(\bar{t}(n))(\gamma_{n+1}-\gamma_n)
\right\}\not=1 
\ee  
on the set of all Fermat's rationals expresses the new $\zeta$-equivalent of the Fermat-Wiles theorem. 
\end{mydef14} 

\section{Influence of the Riemann hypothesis on the $\zeta$-functional (\ref{4.2}) and on the corresponding $\zeta$-equivalent of the Fermat-Wiles theorem}  

\subsection{} 

Let us mention that we have proved, on Riemann hypothesis, the following formula, see \cite{2}, (1): 
\be \label{5.1} 
\sum_{\gamma}\frac{1}{(t-\gamma)^2}=-\frac{{\rm d}}{{\rm d}t}\left\{\frac{Z'(t)}{Z(t)}\right\}+\mcal{O}\left(\frac{1}{t}\right)
\ee 
for every $t\not=\gamma$ and $t\to\infty$. This formula implies the corollary\footnote{See \cite{2}, p. 34.}: The function 
\be \label{5.2} 
\frac{Z'(t)}{Z(t)},\ t\in (\gamma',\gamma'')
\ee 
is decreasing. And further: 
\begin{mydef83}
On Riemann hypothesis the function 
\be \label{5.3} 
|Z(t)|,\ t\in (\gamma',\gamma'') 
\ee  
attains one maximum only. 
\end{mydef83} 

Next, since 
\be \label{5.4} 
\frac{{\rm d}}{{\rm d}t}Z^2(t)=2Z(t)Z'(t),\ t\in (\gamma',\gamma''), 
\ee 
the following property holds true. 

\begin{mydef84}
On the Riemann hypothesis the function 
\be \label{5.5} 
Z^2(t),\ t\in (\gamma',\gamma'')
\ee 
attains one maximum only. 
\end{mydef84} 

Consequently, in the context of the present paper, it holds true: 

\begin{mydef85}
The Riemann hypothesis implies the following behaviour of the function 
\be \label{5.6}  
Z^2(t), \ t\in (\gamma',\gamma''): 
\ee  
\begin{itemize}
	\item[(a)] it has only one point of the maximum 
	\be \label{5.7} 
	\max_{\gamma_n\leq t\leq\gamma_{n+1}}\{Z^2(t)\}=Z^2(t_0(n)),\ \gamma_n<t_0(n)<\gamma_{n+1}, 
	\ee  
	where\footnote{See (\ref{1.15}).} 
	\bdis 
	t_0(n):\ Z'(t_0(n))=0, 
	\edis  
	\item[(b)] and, of course, the function 
	\be \label{5.8} 
	Z^2(t),\ t\in(\gamma_n,t_0) 
	\ee 
	is increasing, and the function 
	\be \label{5.9} 
	Z^2(t),\ t\in(t_0,\gamma_{n+1}) 
	\ee
	is decreasing. 
\end{itemize} 
\end{mydef85} 

\subsection{} 

Let us remind the Bonnet mean-value theorem for the integral 
\be \label{5.10} 
\int_a^b f(x)g(x){\rm d}x . 
\ee  
If $f(x),\ x\in[a,b]$ is decreasing and positive function and $g(x)$ is integrable, then 
\be \label{5.11} 
\int_a^bf(x)g(x){\rm d}x=f(a)\int_a^\xi g(x){\rm d}x,\ \xi\in (a,b), 
\ee  
and if $f(x)$ is incresing, then 
\be \label{5.12} 
\int_a^bf(x)g(x){\rm d}x=f(b)\int_\xi^bg(x){\rm d}x,\ \xi\in (a,b). 
\ee 
Now, by Riemann hypothesis and by Property 5: 
\be \label{5.13} 
\begin{split} 
& \int_{\gamma_n}^{t_0}Z^2(t){\rm d}t=\int_{\gamma_n}^{t_0}Z^2(t)\times 1{\rm d}t=Z^2(t_0)\int_{\xi_1}^{t_0}1{\rm d}t= \\ 
& Z^2(t_0)(t_0-\xi_1),\ \xi_1\in(\gamma_n,t_0), 
\end{split}  
\ee  
and 
\be \label{5.14} 
\begin{split} 
	& \int_{t_0}^{\gamma_{n+1}}Z^2(t){\rm d}t=Z^2(t_0)(\xi_2-t_0),\ \xi_2\in(t_0,\gamma_{n+1}). 
\end{split}  
\ee  
Consequently, we have the following result. 

\begin{mydef53}
On Riemann hypothesis it is true that 
\be \label{5.15} 
\int_{\gamma_n}^{\gamma_{n+1}}Z^2(t){\rm d}t=Z^2(t_0(n))[\xi_2(n)-\xi_1(n)], 
\ee  
where 
\be \label{5.16} 
\gamma_n<\xi_1(n)<t_0(n)<\xi_2(n)<\gamma_{n+1},\ n=1,\dots,N. 
\ee 
\end{mydef53} 

\begin{remark}
The values 
\be \label{5.17} 
\begin{split}
& \xi_1=t_0-\frac{1}{Z^2(t_0)}\int_{\gamma_n}^{t_0}Z^2(t){\rm d}t, \\ 
& \xi_2=t_0+\frac{1}{Z^2(t_0)}\int_{t_0}^{\gamma_{n+1}}Z^2(t){\rm d}t
\end{split}
\ee 
represent the Bonnet points of the integral 
\be \label{5.18} 
\int_{\gamma_n}^{\gamma_{n+1}}Z^2(t){\rm d}t. 
\ee 
\end{remark} 

Next, we obtain the following result from (\ref{3.12}) by (\ref{5.15}). 

\begin{mydef15}
On Riemann hypothesis it is true that 
\be \label{5.19} 
\begin{split}
& \sum_{T<\gamma_n<\overset{1}{T}(T)}Z^2(t_0(n))[\xi_2(n)-\xi_1(n)]=(1-c)T+\mcal{O}(\sqrt{T}),\ T\to\infty, 
\end{split} 
\ee 
where 
\bdis 
\gamma_n<\xi_1(n)<t_0(n)<\xi_2(n)<\gamma_{n+1},\ n=1,\dots,N(T). 
\edis  
\end{mydef15} 

\begin{remark}
The formula (\ref{5.19}) represents, on Riemann hypothesis, the third $\zeta$ sum-variant of our almost linear formula (\ref{3.1}). 
\end{remark} 

\subsection{} 

Now, if we use the substitution (\ref{4.1}) in the formula (\ref{5.19}), then we obtain the next result. 

\begin{mydef16}
The Riemann hypothesis implies the following $\zeta$-sum functional 
\be \label{5.20} 
\lim_{\tau\to\infty}\frac{1}{\tau}\left\{
\sum_{t_0(n)>\frac{x}{1-c}\tau}^{t_0(n)<[\frac{x}{1-c}\tau]^1}Z^2(t_0(n))[\xi_2(n)-\xi_1(n)]
\right\}=x
\ee 
for every fixed $x>0$, where 
\bdis 
\gamma_n<\xi_1(n)<t_0(n)<\xi_2(n)<\gamma_{n+1},\ n=1,\dots,N(T). 
\edis 
\end{mydef16} 

Consequently, in the special case of Fermat's rationals, we have the following result, comp. (\ref{4.2}) -- (\ref{4.4}). 

\begin{mydef17}
On Riemann hypothesis it is true: the $\zeta$-condition 
\be \label{5.21} 
\lim_{\tau\to\infty}\frac{1}{\tau}\left\{
\sum_{t_0(n)>\FRm\frac{\tau}{1-c}}^{t_0(n)<[\FRm\frac{\tau}{1-c}]^1}Z^2(t_0(n))[\xi_2(n)-\xi_1(n)]
\right\}\not=1 
\ee 
on the set of all Fermat's rationals expresses the new $\zeta$-equivalent of the Fermat-Wiles theorem. 
\end{mydef17} 

Now, let the symbol $\mcal{B}(n)$ denote the rectangle 
\be \label{5.22} 
(\xi_1(n),\xi_2(n))\times (0,Z^2(t_0(n))), 
\ee  
where\footnote{See (\ref{3.14}), (\ref{5.14}).} 
\be \label{5.23} 
|\mcal{B}(n)|=Z^2(t_0(n))[\xi_2(n)-\xi_1(n)]=|\mcal{A}_0(n)|. 
\ee 

\begin{remark}
We have substituted, by Riemann hypothesis, the set\footnote{See (\ref{3.13}).}
\bdis 
\bigcup_{n=1}^N\mcal{A}_0(n) 
\edis  
of the curvilinear signals $A_0(n)$ in (\ref{4.4}) by the disconnected set\footnote{See (\ref{5.20}), (\ref{5.22}).} 
\be \label{5.24} 
\bigcup_{n=1}^N\mcal{B}(n) 
\ee 
of isolated rectangular pulses (\ref{5.19}). 
\end{remark} 

\begin{remark}
Of course, the transformation 
\be \label{5.25} 
\mcal{A}_0(n)\to \mcal{B}(n) 
\ee 
preserves the measure together with the transformation\footnote{See (\ref{3.18}).} 
\bdis 
\mcal{A}_0(n)\to\mcal{A}_1(n). 
\edis 
\end{remark} 

\section{Riemann hypothesis and some basic arithmetical and geometrical properties connected with the signals $\mcal{A}_1(n)$ and $\mcal{B}(n)$} 

\subsection{} 

Let us remind two formulae (\ref{3.14}) and (\ref{5.15}): 
\bdis 
\int_{\gamma_n}^{\gamma_{n+1}}Z^2(t){\rm d}t=Z^2(\bar{t}(n))(\gamma_{n+1}-\gamma_n),\ \bar{t}(n)\in(\gamma_n,\gamma_{n+1}),
\edis 
and 
\bdis 
\int_{\gamma_n}^{\gamma_{n+1}}Z^2(t){\rm d}t=Z^2(t_0(n))[\xi_2(n)-\xi_1(n)], 
\edis 
where 
\be \label{6.1}  
\gamma_n<\xi_1(n)<t_0(n)<\xi_2(n)<\gamma_{n+1},\ n=1,\dots,N(T). 
\ee 
These two formulae imply the following result. 

\begin{mydef86}
On Riemann hypothesis it is true that 
\be \label{6.2} 
\frac{Z^2(t_0(n))}{Z^2(\bar{t}(n))}:\frac{\gamma_{n+1}-\gamma_n}{\xi_2(n)-\xi_1(n)}=1
\ee 
for every sifficiently big $T>0$ and every $n=1,\dots,N(T)$, comp. (\ref{3.6}). 
\end{mydef86} 

\begin{remark}
The formula (\ref{6.2}) represents, on Riemann hypothesis, the global conservation law controlling behaviour of four characteristic values: measures of the segments 
\bdis 
[\gamma_n,\gamma_{n+1}],\ [\xi_1(n),\xi_2(n)]
\edis 
together with the values 
\bdis 
Z^2(t_0(n)),\ Z^2(\bar{t}(n)), 
\edis  
i.e. the maximum value and the usual mean value, respectively. 
\end{remark} 

Further, it follows from (\ref{6.2}), for example: 

\begin{mydef87}
On Riemann hypothesis we have the following local proportion 
\be \label{6.3} 
\frac{Z^2(t_0(n))}{\gamma_{n+1}-\gamma_n}=\frac{Z^2(\bar{t}(n))}{\xi_2(n)-\xi_1(n)}=K(n,T), 
\ee  
for every fixed admissible $T$ and $n$. 
\end{mydef87}

I would like to thank Michal Demetrian for his moral support of my study of Jacob's ladders.

\end{document}